\documentclass[leqno,12pt]{amsart}
\usepackage{eepic,epic}
\usepackage{amsmath}

\theoremstyle{plain}
\newtheorem{lemma}{Lemma}[section]
\newtheorem{theorem}[lemma]{Theorem}

\newtheorem*{corollary}{Corollary}

\theoremstyle{definition}

\theoremstyle{remark}
\newtheorem{remark}{Remark}

\def\dim{{\operatorname{dim}}}

\def\exp{{\operatorname{exp}}}

\def\Z{{\mathbb Z}}
\def\Q{{\mathbb Q}}
\def\C{{\mathbb C}}

\def\H{{\mathcal H}}

\def\O{{\mathcal O}}

\def\a{{\mathfrak a}}
\def\pp{{\mathfrak p}}
\begin{document}

\title[Eta-product $\eta(7\tau)^7/\eta(\tau)$]
{
eta-product $\eta(7\tau)^7/\eta(\tau)$}
\author{Kyoji Saito }
\address{Research Institute for Mathematical Sciences, Kyoto University, 
Kyoto 606-8502, Japan.}

\maketitle

\vspace{-0.6cm}
\begin{abstract}
Let $L_{\Phi_7}(s)$ be the Dirichlet series
 associated to the
eta-product $\eta(7\tau)^7/\eta(\tau)\!\in\! M_3(\Gamma_0(7),\varepsilon)$  
(here $\varepsilon(n)\!:=\!\big(\frac{n}{7}\big)\!=\!\big(\frac{-7}{n}\big)$ 
is the Dirichlet character defined by the residue symbol).
We show that $L_{\Phi_7}(s)$ decomposes into the difference of two 
$L$-functions:
\[
L_{\Phi_7}(s)=\frac{1}{8}\big(L(s,\varepsilon)L(s-2,1)-L(s-1,\xi)\big),
\]
where  
i) $L(s,\varepsilon)$ and $L(s,1)$ are Dirichlet $L$-functions 
for the characters $\varepsilon$ and 1 modulo 7, respectively, 
and ii) $L(s,\xi)$ is the $L$-function for a Hecke character $\xi$ 
of the imaginary quadratic field $\Q(\sqrt{-7})$.

This expression of $L_{\Phi_7}(s)$ gives a new proof of the non-negativity of 
the Fourier coefficients of the product $\eta(7\tau)^7/\eta(\tau)$, 
conjectured in [S3] and proven by Ibukiyama [I].
We also prove the uniqueness of the above decomposition of 
$L_{\Phi_7}(s)$ in a suitable sense.
\end{abstract}

\vspace{-0.2cm}
\section{introduction}

Let 
$\eta(\tau)\!=\! 
q^{\frac{1}{24}}\prod_{n=1}^{\infty} (1\!-\!q^{\eta}), 
\ q\! =\! \exp (2 \pi \sqrt{-1} \tau)$ be the Dedekind eta-function
(e.g.\ [R]). 
A product $\Pi_{i\in I}\eta(i\tau)^{e(i)}$, where 
$I$ is a finite set of positive integers and $e:I\to \Z$ is any map,
 is called an eta-product.
The eta-product can be developed in a Laurent series in powers of $q$,
whose coefficients are called the {\it Fourier coefficients}.

Ibukiyama [I] has shown the following result, which answers to a part 
of a conjecture given by the author [S3] (see the next paragraph).

\begin{theorem}
Let $p$ be a rational prime number. Then the Fourier coefficients of the 
eta product $\eta_{\Phi_p}:=\eta(p\tau)^p/\eta(\tau)$ are non-negative.
\end{theorem}

\noindent
The proof in [I] is given by expressing the eta-product as a 
difference of two generating functions of two arithmetically constructed 
lattices.

 More in general than the theorem, 
for any positive integer $h$ which may not be prime, 
 we have the following non-negativity conjecture.

\medskip
\noindent
{\footnotesize
{\it Conjecture} ([S3]).
Define a sequence $\Phi_h(\lambda)$ ($h\in \Z_{>0}$) of cyclotomic polynomials 
by the recursive relation: 
$
 \frac{(1-\lambda^h)^h}{1-\lambda} = \prod_{d|h}\Phi_d(\lambda^{h/d}).
$
Explicitly, 
$
\Phi_h(\lambda)= \frac{(1-\lambda^h)^{\phi(h)}}{\prod_{d|h}(1-\lambda^d)^{\mu(d)}}
$ 
where $\phi$ and $\mu$ are the Euler function and the M\"obius function. 
Then the Fourier coefficients of 
the eta-product
$
\eta_{\Phi_h}(\tau):=\frac{\eta(h\tau)^{\phi(h)}}
{ \prod_{d|h}\eta(d\tau)^{\mu(d)}}
$
are non-negative integers.
}

\medskip
\noindent
This was proven for
 $h\!=\!2,3,4,5,6$ [S1,2,3] by a use of the 
Dirichlet series $L_{\Phi_h}(s)$ associated to 
the eta-products $\eta_{\Phi_h}$.	
Precisely, we show that 
$L_{\Phi_h}(s)$ admits either an Euler product for $h=2,3,5$ or  
a decomposition into a difference of two 
Euler products for $h=4,6$, and that  
these expressions lead to a direct proof of the positivity 
of the coefficients.


In the present note, we prove in section 2 that the Dirichlet series 
$L_{\Phi_7}(s)$ decomposes into a difference of two $L$-functions,
which admit Euler products, as stated in Abstract.
In section 3, we show that this expression implies the 
non-negativity of the Dirichlet coefficients of $L_{\Phi_7}(s)$. 
In section 4, we prove a general lemma on the uniqueness of 
the decomposition of Dirichlet series into a difference of two 
Euler products, and apply it to $L_{\Phi_7}(s)$ (and also to 
$L_{\Phi_4}(s)$ and $L_{\Phi_6}(s)$). 
Finally, we remark in section 5 that such difference decomposition of  
$L_{\Phi_p}(s)$ for the prime $p\ge11$ does not exist.
If $h$ is a composite number, we  do not know when 
$L_{\Phi_h}(s)$ admits such a difference decomposition.

In [S2,Conjecture 13.5], we give a wide class of eta-products whose 
Fourier coefficients are conjecturally non-negative and are of interest.

\vspace{-0.1cm}
\section{Hecke $L$-function $L(s,\xi)$ for a character $\xi$ on $\Q(\sqrt{-7})$}

We recall Hecke's $L$-function for a character $\xi$ on the imaginary 
quadratic field 
$\Q(\sqrt{-7})$, and, then, decompose $L_{\Phi_7}(s)$ by a use of it. 
For a back ground on analytic number theory, one is referred to
[M] and [R].

Since the class number of $\Q(\sqrt{-7})$ is equal to 1, 
we can introduce the Hecke character $\xi$ for the non-zero ideals of 
$K:=\Q(\sqrt{-7})$ by 
\begin{equation}
\xi((a)):=\big(\frac{a}{|a|}\big)^2 \qquad (a\in K\setminus\{0\}).
\end{equation}
Then, the $L$-function for $\xi$ is defined by the following Dirichlet series, 
which, as a result of definition, has the Euler product:
\begin{equation}
 L(s,\xi):=\sum_{\a\ \subset\O_K} \xi(\a)N_K(\a)^{-s}=\prod_{\pp\ :\ {\rm prime}}(1-\xi(\pp)N_K(\pp)^{-s})^{-1}.
\vspace{-0.1cm}
\end{equation} 
Here, $\a$ (resp. $\pp$) runs over all non-zero integral (resp. prime) 
ideals of $\O_K$, and $N_K(\a)$ is the absolute 
norm of $\a$ (i.e.\ $N_K(\a)=|\O_K/\a|$). 

\medskip
The first main result of the present note is the following.
\begin{lemma} The Dirichlet series $L_{\Phi_7}(s)$ 
associated to the eta-product $\eta(7\tau)^7/\eta(\tau)$ 
decomposes into a difference of two $L$-functions as follows:
\begin{equation}
L_{\Phi_7}(s)=\frac{1}{8}\big(L(s,\varepsilon)L(s-2,1)-L(s-1,\xi)\big),
\end{equation}
where we recall that 
$\varepsilon\!:\!=\!\big(\frac{*}{7}\big)\!=\!\big(\frac{-7}{*}\big)$
 is the residue symbol modulo 7.
\end{lemma}
\begin{proof}

The $L$-function  $L(s-1,\xi)$ is associated to a Fourier series 
\begin{equation}
f(\tau):=\sum_{\a}\xi(\a)N_K(\a)e^{2\pi\sqrt{-1}N_K(\a)\tau}.
\end{equation}
According to Hecke [H1][H2], $f(\tau)$ is an automorphic form 
belonging in $S_3(\Gamma_0(7),\varepsilon)$ (see [M, Th.4.8.2]).  
Similarly, $L(s,\varepsilon)L(s-2,1)$ and $L(s-2,\varepsilon)L(s,1)$ 
are associated to Eisenstein series, say $E(\tau)$ and $E'(\tau)$, 
in $M_3(\Gamma_0(7),\varepsilon)$. 
Since $\Gamma_0(7)\backslash\H$ has two cusps and $\dim_\C S_3(\Gamma_0(7),\varepsilon)\!=\!1$, $M_3(\Gamma_0(7),\varepsilon)$ is spanned by $E,E'$ and $f$.
To show the equality:
$\eta_{\Phi_7}(\tau)=\frac{1}{8}\big(E(\tau)-f(\tau)\big)$,
it suffices to show that 
$n$th Fourier coefficients $c(n)$ of $\eta_{\Phi_7}(\tau)$ 
coincide with $n$th Dirichlet coefficients of 
$\frac{1}{8}(L(s,\varepsilon)L(s\!-\!2,1)\!-\!L(s\!-\!1,\xi))$ 
for $1\!\le\!n\!\le\!3$. 
We give an explicit integral description 
(which we shall use in the next section) of the coefficients of $L(s-1,\xi)$.
For the end,
{\it we factorize $L(s-1,\xi)$ w.r.t.\ 
rational primes $p,q$ in $\Z_{>0}$:
\begin{equation}
L(s-1,\xi):=\frac{1}{1+7^{-s+1}}\cdot 
\prod_{\varepsilon(q)=-1}\frac{1}{1-q^{-2s+2}} \cdot
\prod_{\varepsilon(p)=1}\frac{1}{P_p(p^{-s})},
\end{equation}
where $P_p(\lambda)\! \in\!\Z[\lambda]$ 
for a prime $p$ with $\varepsilon(p)\!=\!1$ is defined 
in next} (6). 

\medskip
\noindent
{\it Proof.}
Recall a well-known fact (e.g.\ [T]) on the prime ideals in $\Q(\sqrt{-7})$:

i) $(q)$ is a prime ideal for any rational prime $q$ with $\varepsilon(q)=-1$, 

ii)
$
p=x_p^2+ 7\cdot y_p^2= (x_p+y_p\sqrt{-7})(x_p-y_p\sqrt{-7})
$ \quad 
($(x_p,y_p)\in \Z_{>0}^2$)
for any odd rational prime number $p$ with $\varepsilon(p)=1$,

iii) $2 =\frac{7\cdot1 + 1}{4}= \frac{1+\sqrt{-7}}{2}\cdot \frac{1-\sqrt{-7}}{2}$ and $7=-(\sqrt{-7})^2$. 
\quad 

Put
$
\pi_2:=\frac{1+\sqrt{-7}}{2} \text{\ and } \pi_p:=x_p+y_p\sqrt{-7}
$
for an odd rational prime number $p$ with $\varepsilon(p)=1$ 
and, define the quadratic polynomials 
\begin{equation}
\begin{array}{lllll}
P_2(X)\!&:=&\!(1-\pi_2^2X)(1-\overline{\pi}_2^2X)\!&=&\!1 + 3X+2^2 X^2 \text{  and}\\
P_p(X)\!&:=&\!(1-\pi_p^2X)(1-\overline{\pi}_p^2X)\!&=&\!1 - 2(x_p^2-7y_p^2)X+p^2 X^2.\!\!
\end{array}\!\!
\end{equation}
Then (5) follows from the Euler product in (2) and  

i) \ \  $\xi((\pi_p))=\pi_p^2/p$ and $N_K((\pi_p))=p$ for $\varepsilon(p)=1$, 

ii) \ $\xi((q))=1$ and $N_K((q))=q^2$ for $\varepsilon(q)=-1$,
 
iii) $\xi((\sqrt{-7}))=-1$ and $N_K((\sqrt{-7}))=7$ 
. \qquad\qquad \qquad\qquad \qquad $\Box$

\medskip
Put $L(s,\varepsilon)L(s-2,1)=\sum_{n=1}^\infty a(n)n^{-s}$ and 
$L(s-1,\xi)=\sum_{n=1}^\infty b(n)n^{-s}$,
and we give explicite expressions of the coefficients $a(n)$ and
$b(n)$. Let

\centerline{
$n=7^k\prod_{i\in I}p_i^{l_i}\prod_{j\in J}q_j^{m_j}$
}

\noindent
be the prime decompostion of $n\!\in\!\Z_{>0}$ where $\{p_i\mid i\!\in\! I\}$ and 
$\{q_j\mid j\!\in\! J\}$ are finite set of distinct prime numbers 
with $\varepsilon(p_i)=1$ and $\varepsilon(q_j)=-1$. 

Then, by a use of (5) together with (6), one obtains the formulae:
\begin{eqnarray}
a(n)&=&7^{2k}\prod_{i\in I}\frac{p^{2(l_i+1)}-1}{p_i^2-1}
\prod_{j\in J}\frac{q_j^{2(m_j+1)}-(-1)^{m_j+1}}{q_j^2+1} \\
b(n)&=&(-7)^{k}\prod_{i\in I}\big(\sum_{t=0}^{l_i}\pi_{p_i}^{2t}\overline{\pi}_{p_i}^{2(l_i-t)}\big)
\prod_{j\in J}\frac{1-(-1)^{m_j+1}}{2}q_j^{m_j} 
\end{eqnarray}

Finally, we give the Fourier expansion of $\eta_{\Phi_7}$ up to degree 50. 
{\scriptsize
\[\begin{array}{rll}
\!\!\!\eta_{\Phi_7}\!\!&\!\!=\!\!&\!q^2+q^3+2q^4+3q^5+5q^6+7q^7+11q^8+8q^9+15q^{10}+16q^{11}+21q^{12}+21q^{13}\\
&&\!\!\!\!\!\!\!+28q^{14}+24q^{15}+44q^{16}+36q^{17}+49q^{18}+45q^{19}+63q^{20}+49q^{21}+74q^{22}+64q^{23}\\
&&\!\!\!\!\!\!\!+85q^{24}+72q^{25}+105q^{26}+82q^{27}+133q^{28}+112q^{29}
+120q^{30}+120q^{31}+165q^{32}\\
&&\!\!\!\!\!\!\!+122q^{33}+180q^{34}+147q^{35}+186q^{36}+176q^{37}+225q^{38}
+168q^{39}+255q^{40}+21q^{41}\\
&&\!\!\!\!\!\!\!+245q^{42}+224q^{43}+324q^{44}+219q^{45}+338q^{46}
+276q^{47}+341q^{48}+294q^{49}+385q^{50}+\cdots
\end{array}
\]
}
By inspection, we check the equality $c(n)=\frac{1}{8}(a(n)-b(n))$ for $n$ 
with $1\le n\le 3$. This completes a proof of Lemma 2.1.
\end{proof}
\begin{remark}
As we see in  the above proof, once one guess a correct 
formula (3), then its proof 
is straight forward. However, we do not know yet what is a 
``correct formula'' for $L_{\Phi_h}(s)$ for $h>7$ (see \S 5).
\end{remark}

\section{Positivity of Fourier coefficients of $\eta(7\tau)^7/\eta(\tau)$}
As an immediate consequence of Lemma 2.1.\  together with the 
explicit formulae (6) and (7), we obtain the following positivity.

\begin{corollary} All Fourier coefficients 
of $\eta(7\tau)^7/\eta(\tau)$ are positive.
\end{corollary}
\begin{proof}
Lemma 2.1.\ says $c(n)=\frac{1}{8}(a(n)-b(n))$ for all $n\in\Z_{\ge1}$. 
To show $a(n)>b(n)$ for all $n\in\Z_{\ge1}$, it is sufficient to show 
$a(p^k)>|b(p^k)|$ for any primary number $p^k$ (i.e. $p$ is a prime 
number and $k\in\Z_{>0}$) because of the multiplicativity of $a(n)$ and $b(n)$.
We separate cases:

\noindent
Case  $p=7$.   \
$a(7^k)=7^{2k}>7^k=|b(7^k)|$.

\noindent
Case  $\varepsilon(p)=1$. \  
$a(p^k)>p^{2k}\ge (k+1)p^k=\sum_{i=0}^k|\pi_p^{2i}\overline{\pi}_p^{2(k-i)}|
\ge |b(p^k)|$.

\noindent
Case $\varepsilon(q)\!=\!-1$.  
$a(q^k)\!-\!|b(q^k)|\ge\frac{q^{2(k+1)}-1}{q^2+1}-q^{k}\!=\!\frac{(q^{k+2}-1)(q^{k}-1)-2}{q^2+1}\!>\!0$.
\end{proof}

\section{Uniqueness of decomposition of Dirichlet series}

We show the second main result of the present note: 

\noindent
Under a mild assumption on 
a Dirichlet series $L(s)=\sum_{n\in\Z_{\ge1}}c(n) n^{-s}$, we show 
{\it the uniqueness of the decomposition of 
$L(s)$ into the form: 
\begin{equation}
L(s)=aM(s)+bN(s)
\end{equation}
where  $M(s)$ and $N(s)$ are Dirichlet series which 
admit Euler product and $a,b$ are constants.} 
For our applications, we assume that $c(1)=0$ so that 
one automatically has $a+b=0$ (since the first Dirichlet 
coefficients of $M(s)$ and $N(s)$ are automatically equal to 1).

\begin{lemma} 
Let $L(s)=\sum_{n\in\Z_{\ge1}}c(n) n^{-s}$ be a Dirichlet series
such that {\rm i)} $c(1)=0$ and {\rm ii)}  there are five relatively prime 
integers $l,m,n,u,v\in\Z_{\ge1}$ such that $c(l)c(m)c(n)c(u)c(v)\not=0$. 
If there exists a decomposition {\rm (9)}, where $M(s)$ and $N(s)$ 
are Dirichlet series 
having Euler products, then it is unique
up to the transposition of $M(s)$ and $N(s)$.
\end{lemma}
\begin{proof}
Put $M(s)=\sum_{n\in\Z_{\ge1}}a(n) n^{-s}$,
$N(s)=\sum_{n\in\Z_{\ge1}}b(n) n^{-s}$ and $c:=a=-b$ 
so that one has the relation among the Dirichlet coefficients:
\begin{equation}
\qquad c(n)=c( a(n)-b(n)) \qquad (n\in\Z_{\ge1}). 
\end{equation}
Clearly $c\not=0$, else $L(s)=0$ contradicting to the assumption on $L(s)$.

We first remark that one sees from (10) that if $c(n)=c(m)=0$ 
for relatively prime positive integers $n$ and $m$ then $c(nm)\!=\!0$. 
Consequently, if $c(n)\!\not=\!0$, then there exists a primary 
factor $p^k$ of $n$ 
(i.e.\ $p$ is a prime number and $k$ is a positive integer s.t.\ $p^k|n$) 
such that $c(p^k)\not=0$.

Suppose there exist another decomposition
$L(s)=c'(M'(s)-N'(s))$. 
Using Dirichlet coefficients $a'(n),b'(n)$ of 
$M'(s),N'(s)$, this means
\begin{equation}
\qquad c(n)=c'(a'(n)-b'(n))  \qquad (n\in \Z_{\ge1})
\end{equation}

Let $n,m\in\Z_{\ge1}$ be relatively prime to each other, 
then the multiplicativities 
of the Dirichlet coefficients $a,b,a'$ and $b'$ implies 
\[
c(mn)=c(a(n)a(m)-b(n)b(m))=c'(a'(n)a'(m)-b'(n)b'(m))
\]
Substituting $b(n)=a(n)\!-\!c(n)/c$, $b'(n)=a'(n)\!-\!c(n)/c'$ and  
$b(m)$ $=a(m)\!-\!c(m)/c$, $b'(m)=a'(m)\!-\!c(m)/c'$ in *),
we obtain 
\[
c(n)(a(m)\!-\!a'(m))+c(m)(a(n)\!-\!a'(n))=(\frac{1}{c}\!-\!\frac{1}{c'})c(n)c(m).
\leqno{E(m,n):}
\]
Let $k,m,n\!\in\! \Z_{\ge1}$ be relatively prime to each other and
$c(m)c(n)\!\not=\!0$, 

\noindent
then
$(c(k)E(m,n)\!-\!c(m)E(n,k)\!-\!c(n)E(k,m))/c(m)c(n)$ is the equality
\[
a(k)-a'(k)=\frac{1}{2}(\frac{1}{c}-\frac{1}{c'})c(k).
\leqno{*}
\]
This, together with (10) and (11), can be rewritten as the linear relations 
among $a(k),b(k)$ and $a'(k),b'(k)$ for all $k$ prime to $mn$:
\[
a'(k)=(1-\lambda) a(k)+\lambda b(k) \quad \text{and}\quad
b'(k)=\lambda a(k)+(1-\lambda) b(k),
\]
where $\lambda:=\frac{c}{2}(\frac{1}{c}-\frac{1}{c'})$ so that 
$\lambda=0$ or $1$ if and only if $c=c'$ or $c=-c'$, respectively. 
Summing two relations, we also obtain the relation:
\[
a(k)+b(k)=a'(k)+b'(k).
\leqno{**}
\]
If $c=c'$ (i.e.\ $\lambda=0$), then the proof of Lemma 4.1.\ is already achieved 
as follows: by substituting $c=c'$ in $*$ and using $**$, one has 
\[
a(k)=a'(k) \text{\quad and \quad} b(k)=b'(k)
\leqno{***}
\] 
for any $k\in\Z_{\ge1}$ prime to $m,n$. By replacing the role of $m,n$  
by $u,v$, the equalities $***$ hold for any primary numbers $k$. 
The $***$ extends, further, for any positive integers $k$ due to 
the multiplicativity of $a,a',b$ and $b'$. This means
$M(s)=M'(s)$ and $N(s)=N'(s)$.

Suppose $c\not=c'$ (i.e.\ $\lambda\not=0$). 
Then, $*$ means another decomposition:
\[
c(k)=\frac{c}{\lambda}(a(k)-a'(k))
\leqno{(11)'}
\]
for any $k\in\Z_{\ge1}$ prime to $m,n$.  Replacing (11) by (11)', we can 
repeat the previous discussions to induce $*$ and $**$, where we replace
the role of $m$, $n$ by $u$, $v$,
and consider integers $k$ which is prime to $m,n$ and
also to $u,v$. Then, in addition to $*$ and $**$, we obtain:
$*':\ \ 0=a(k)-a(k)=\frac{1-\lambda}{2c}c(k)$ and
$**':\ \ a(k)+b(k)=a(k)+a'(k)$
for all $k$ prime to $m,n,u,v$. 
Taking $k=l$ with $c(l)\not=0$, 
which exists by the assumption of Lemma,
we obtain $\lambda\!=\!1$, i.e.\ $c\!=\!-c'$. 
By the similar argument for the case $c\!=\!c'$, 
we obtain: $***':\ \ a(k)\!=\!b'(k),\ b(k)=a'(k)$ for all $k\in\Z_{\ge1}$ 
and, therefore, $M(s)=N'(s)$ and $N(s)=M'(s)$.
\end{proof}
\begin{corollary}
The Dirichlet series $L_{\Phi_7}(s)$ satisfies the assumptions 
{\rm i)} and {\rm ii)} 
so that the decomposition {\rm (3)} is unique in the sense of Lemma {\rm 4.1}.
\end{corollary}

\begin{remark}
Lemma 4.1.\ can be formulated more precisely according to the \# of  
relatively prime $n$'s with $c(n)\not=0$.
The case \#=5 of Lemma 4.1.\ is the strongest case.
Since the other cases for $\# <5$ are involved but not used in the present note, 
they are omitted.
\end{remark}
\begin{remark} 
There are a few more known Dirichlet series associated to  
eta-products, which decompose as (9) and 
satisfy the assumption of Lemma 4.1, namely, 
$\eta(48\tau)^3/\eta(24\tau)$, 
$\eta_{\Phi_4}(8\tau)=\eta(32\tau)^2\eta(16\tau)/\eta(8\tau)$ and 
$\eta_{\Phi_6}(12\tau)=\eta(72\tau)\eta(36\tau)\eta(24\tau)/\eta(12\tau)$.
They have an origin in a study of elliptic root systems (see [S1]).

\end{remark}

\section{non-decomposability of $L_{\Phi_p}(s)$ for $p\ge11$}
We finally give the following remark, which can be shown trivially.

\medskip
\noindent
{\bf Fact.} {\it The Dirichlet series $L_{\Phi_p}(s)$ associated to a eta-product 
$\eta(p\tau)^p/\eta(\tau)$ for a prime number $p$ with $p\ge11$ does not 
admit a decomposition {\rm (9)}.}

\begin{proof}
Suppose a decomposition (9) exists, i.e.\ there is a Dirichlet series $M(s)$ 
and a constant $c\not=0$ such that $M(s)-\frac{1}{c}L_{\Phi_p}(s)$ is a 
Dirichlet series admitting an Euler product.
Let $c(n)$, $a(n)$ and $b(n)$ be the Dirichlet coefficients of $L_{\Phi_p}(s)$, 
$M(s)$ and $M(s)-\frac{1}{c}L_{\Phi_p}(s)$.  The following fact follows 
from the explicit
expression of the eta product $\eta(p\tau)^p/\eta(\tau)$:

i) \  $c(n)=0$ \ for $1\le n<(p^2-1)/24\ (\ge 5)$, 
 
ii) $c(n)\not=0$ \ for $(p^2-1)/24\le n<(p^2-1)/24+p$.

Thus, we can find an odd integer $m$ such that $1<m<(p^2-1)/24$ and 
$(p^2-1)/24\le 2m<(p^2-1)/24+p$. Then, 
$a(2)a(m)=b(2)b(m)=b(2m)=a(2m)-\frac{1}{c}c(2m)=a(2)a(m)-\frac{1}{c}c(2m)$
should imply $\frac{1}{c}c(2m)=0$. Since $c(2m)\not=0$ (due to ii)), 
one has $\frac{1}{c}=0$ which is impossible.
\end{proof}

\noindent
{\bf Acknowledgement} : The author is grateful to Professor 
Hiroshi Saito for his help to identify $L(s,\xi)$ 
with Hecke $L$-function.

\bigskip

\centerline{
REFERENCES}

\medskip\noindent
[H1] Hecke, Erich, \"Uber die Bestimmung Dirichletscher Reihen durch ihre Funktionalgleichung, Math. Ann.112, (1936) \S7.

\medskip\noindent
[H2] Hecke, Erich, \"Uber Modulfunktionen und die Dirichletschen Reihen mit Eulerscher Produktenwicklung II, Math.Ann. {\bf 114} (1937) \S13.


\medskip\noindent
[I2] Ibukiyama, Tomoyosi, Positivity of Eta Products - a Certain Case of K. Saito's Conjecture, Publ. RIMS, Kyoto Univ. {\bf 41} (2005), 683-693.

\medskip\noindent
[M] Miyake, Toshitsune, Modular Forms, Springer Verlag Berlin 
Heidelberg New York 1989, 2006.


\medskip\noindent
[R] Rademacher, Hans, Topics in analytic Number Theory, Springer-Verlag Berlin Heidelberg New York 1973, 3-540-o5447-2.

\medskip\noindent
[S1] Saito, Kyoji, Extended Affine Root Systems,
V. (Elliptic Eta-Products and Their Dirichlet Series), Centre de Recherches Mathematiques CRM Proceedings and Lectures Notes Vol. {\bf 30},2001,pp185-222.

\medskip\noindent
[S2] Saito, Kyoji, Duality for Regular Systems of Weights, Asian J. Math. Vol.{\bf 2}, No.4, pp.983-1048, 1998.

\medskip\noindent
[S3] Saito, Kyoji, Non-negativity of fourier Coefficients of Eta-product 
(in Japanese), Proceedings of the second Spring conference on automorphic 
formes and related subjects, Careac Hamamatsu, Feb.\ 2003.

\medskip\noindent
[T] Takagi, Teiji, Algebraic Number Thoery (in japanese), Iwanami, .

\end{document}